\documentclass[12pt]{amsart}
\usepackage{graphics}

\begin{document}

\title[Absolutely continuous spectrum: optimal results]
{Absolutely continuous spectrum for one-dimensional 
Schr\"odinger operators 
with slowly decaying potentials: some optimal results}

\author{Michael Christ}
\address{\hskip-\parindent Michael Christ\\ Department of Mathematics \\
University of California \\
Berkeley CA 94720}
\email{mchrist@math.berkeley.edu}
\author{Alexander Kiselev}
\address{\hskip-\parindent
Alexander Kiselev\\
Mathematical Sciences Research Institute \\
1000 Centennial Drive \\
Berkeley CA 94720 }
\email{kiselev@msri.org }


\begin{abstract} 
The absolutely continuous spectrum 
of one-dimensional Schr\"odinger operators
is proved to be stable under perturbation by potentials satisfying
mild decay conditions. In particular, 
the absolutely continuous spectrum of free and periodic 
Schr\"odinger operators is preserved under all 
perturbations $V(x)$ satisfying $|V(x)|\leq C(1+x)^{-\alpha},$ 
$\alpha >\frac{1}{2}.$ This result is optimal in the power 
scale.
More general classes of perturbing
potentials which are not necessarily power decaying are also treated. 
A general criterion for stability of the absolutely 
continuous spectrum 
of one-dimensional Schr\"odinger operators is established. 
In all cases analyzed,
the main term of the asymptotic behavior of 
the generalized eigenfunctions is shown to have
WKB form for almost all energies.
The proofs rely on new maximal function and norm estimates and 
almost everywhere
convergence results for certain multilinear integral operators. 
\end{abstract}

\maketitle

\begin{center}
 \large \bf 1. Introduction and main results.
\end{center}

In this paper, we study the stability of the absolutely continuous
 spectrum of one-dimensional Schr\"odinger operators under perturbations 
by slowly decaying potentials. Suppose that $H_{U}$ is a Schr\"odinger 
operator defined on $L^{2}(0,\infty)$ by the 
differential expression 
\[ H_{U}= -\frac{d^{2}}{dx^{2}} + U(x) \]
and some self-adjoint boundary condition at the origin. We assume that $U$ is
 some bounded function for which $H_{U}$ has absolutely 
continuous spectrum. The presence of the absolutely continuous
spectrum has direct consequences for the physical properties of the 
quantum particle described by the operator $H_{U}$ (see, e.g.\
\cite{ReSi}, \cite{AvS}). If we
perturb this operator by some decaying potential $V(x),$ 
the Weyl criterion implies that the essential spectra of the operators 
$H_{U}$ and $H_{U+V}$ coincide.  
We seek conditions on the rate of decay of $V(x)$ which 
ensure that the absolutely continuous spectrum of the unperturbed operator
$H_{U}$ is also preserved. 
 
This problem has a long history as one of the most natural questions in 
quantum mechanics, and we briefly recall the main results. It has long been 
known that if the perturbation $V(x)$ is absolutely
integrable, then the absolutely
 continuous spectrum of the original operator is preserved. 
Until recently, little more
was known concerning the preservation of the absolutely continuous spectrum of
Schr\"odinger operators
under decaying perturbations in the general situation.

Substantially more information is available in the case when $U(x)=0.$
There has been much work on proving the
absolute continuity of the spectrum for the Schr\"odinger operators with
potentials of slower decay, but satisfying some additional special  
assumptions.
For example, by a result going back to Weidmann \cite{Weid}, if a
potential $V$ may be represented as a sum of a function of
bounded variation and an absolutely integrable function, then the spectrum
of the operator $H_{V}$ on $R^{+}=(0, \infty)$ is purely absolutely
continuous. Many authors developed
a scattering theory for long-range potentials whose derivatives
satisfy certain bounds;
see for example \cite{As}, \cite{Bus}, \cite{Hor}.  These 
results hold in any
dimension and the proofs involve approximating the scattering
trajectories
by solutions of the classical Hamilton-Jacobi equation. The weakest 
conditions on the long-range 
part of the potential under which the
wave operators are known to exist, are given in \cite{Hor}. For potentials
satisfying $|V(x)| \leq C(1+|x|)^{-\frac{1}{2}-\epsilon},$ for instance, 
one can infer the existence of the wave operators if also 
$|D^{\alpha}V(x)| \leq C_{1}
(1+|x|)^{-\frac{3}{2}-\epsilon}$ for every multiindex 
$\alpha$ with $|\alpha|= 1.$

Another class of results describes spectral behavior of  specific 
spherically symmetric (i.e.\ essentially one-dimensional) oscillating 
potentials, the
typical example being $V(x) = x^{-\beta} {\sin x^{\alpha}}$
with $\alpha,$ $\beta$ positive. Such potentials
in general do not satisfy  the derivative bounds
needed for the method of 
the works cited above to be applicable. We mention the papers  \cite{Ben},
\cite{BA}, \cite{HiS}, \cite{Mat} and \cite{White} in which further 
references may be found. The spectrum of the operator $H_{V}$ for
such potentials turns out to be absolutely continuous with perhaps some 
isolated
embedded eigenvalues when $\alpha =1.$ These potentials generalize the 
celebrated
Wigner-von Neumann example \cite{Wig}. Wigner and von Neumann were the
first 
to discover an example 
with Coulomb type decay at infinity, i.e.\ $V(x)=O(\frac{1}{1+|x|}),$
whose spectrum is not purely absolutely continuous and has positive 
eigenvalues embedded in the absolutely continuous spectrum. 
Moreover, Naboko \cite{Nab} and  later Simon
\cite{Sim1} found different constructions which show that for potentials 
decaying more slowly but arbitrarily close to
a Coulomb rate, very striking spectral phenomena arise.
Namely, for every function $C(x)$ tending monotonically to
infinity as $x$ goes to infinity, no matter how slowly, there exists
a potential $V(x)$ satisfying
\[ |V(x)| \leq \frac{C(x)}{1+|x|}, \]
for which the associated Schr\"odinger operator $H_{V}$ has a dense set 
of eigenvalues in $R^{+}.$

A new general class of potentials preserving the absolutely continuous
spectrum of the free Schr\"odinger operator was recently found by one 
of us in 
\cite{Kis}.
 Namely, if the potential $V$ satisfies
 $|V(x)| \leq C(1+|x|)^{-\frac{3}{4}-\epsilon}$ with some $\epsilon > 0,$
 with no additional assumptions, then the whole 
positive semi-axis $(0,\infty)$ is an essential support
of the absolutely continuous part of the spectral measure.\footnote{
$E$ is an essential support of the measure $\mu$
if $\mu(X)>0$ for any $X \subset E$ of positive Lebesgue measure.}
Of course, as the examples of
Naboko and Simon show, rich embedded singular spectrum may occur; however 
it is indeed embedded in the sense that there is an underlying absolutely 
continuous spectrum. One can describe the set where the singular part of
the spectral measure might be supported in $R^{+}$ rather explicitly in 
terms of the properties of the Fourier transform of 
$x^{\frac{1}{4}}V(x)$ \cite{Kis}.  

The result of \cite{Kis} was further improved in \cite{Kis1}, where a 
general criterion was established which implies 
the stability of the absolutely 
continuous spectrum of the operator $H_{U}$ under all perturbations $V(x)$
satisfying $|V(x)| \leq C(1+x)^{-\frac{2}{3}-\epsilon}$,
under the auxiliary hypothesis that a
certain operator, constructed from the generalized eigenfunctions of $H_{U},$ 
is bounded on $L^2(R)$.
In particular, it was shown that the absolutely 
continuous spectra of free and periodic one-dimensional 
Schr\"odinger operators are stable under all perturbations 
by potentials satisfying
 $|V(x)|\leq C(1+x)^{-\frac{2}{3}-\epsilon}.$ Later, this result for the
 case $U=0$ was also proved by S.~Molchanov by a different method 
\cite{Mol1}.

On the other hand, there exists work on random potentials by
 Kotani and Ushiroya \cite{Kot}  which  provides 
a bound for the
 best possible result that one
can hope to prove. The results of \cite{Kot} 
imply that (in the case $U=0$) there
exist potentials $V(x)$ satisfying  $V(x) \leq C(1+x)^{-\frac{1}{2}},$ 
for which $H_V$ has purely singular spectrum on $R^{+}.$ 
We remark that by combining  methods  used in recent works \cite{KLS}
 and \cite{KRS}, one can show that the $(1+x)^{-\frac{1}{2}}$ 
rate of decay is also critical for perturbations of the periodic
 Schr\"odinger operators.   

In this paper we establish results on the preservation of absolutely
continuous spectrum for power decaying potentials in one dimension
which are of an optimal nature.
In particular, we prove \\

\noindent \bf Theorem 1.1. \it Suppose that 
there exist $\epsilon>0$ and $p\leq 2$ so that $V(x)(1+x)^{\epsilon} 
\in L^{p}(R^{+}).$
Then the whole positive semi-axis $R^{+}=(0,\infty)$ is an essential
 support of the absolutely continuous part of the spectral measure of the
 operator $H_{V}.$ 
Moreover, for almost every $\lambda \in R^{+},$ there exist solutions 
$\phi_{\lambda}(x)$ and $\overline{\phi}_{\lambda}(x)$ of the generalized 
eigenfunction equation 
\[ -y'' + V(x)y=\lambda y \]
with asymptotic behavior of pure WKB form in the main term:
\[ \phi_{\lambda}(x)= \exp \left(i\sqrt{\lambda}x - 
\frac{i}{2\sqrt{\lambda}}\int\limits_{0}^{x}V(t)\,dt\right)(1+o(1)). \]

\noindent \bf Corollary. \it 
If $V(x) = O(1+x)^{-r}$ for some $r>1/2$ then an essential support of
the absolutely continuous part of the spectral measure
 of $H_V$ on $L^2(R^+)$
equals $[0,\infty)$.
\\

\rm
We have a similar result for perturbations of periodic 
Schr\"odinger operators. Let $U$ be continuous and periodic, 
let $S=\cup_{n=0} ^{\infty}(a_{n},b_{n})$ be
the band spectrum of the unperturbed operator with potential $U$, 
and let $\theta(x, \lambda)$ be the Bloch functions for that operator. 
\\

\noindent \bf Theorem 1.2. \it If the potential $V(x)$ is as in Theorem
 1.1, then the set $S$ is an essential support of the absolutely 
continuous part of the spectral measure of the operator $H_{U+V}.$
 For a.e.\ $\lambda \in S,$ there exist solutions
 $\psi_{\lambda}(x),$ $\overline{\psi}_{\lambda}(x)$ of the equation 
\begin{equation}
-y''+(U(x)+V(x))y = \lambda y
\end{equation}
with the asymptotic behavior
\begin{equation}
 \psi_{\lambda}(x)= \theta(x,\lambda)\exp \left( 
\frac{i}{2\Im(\theta\overline{\theta}')}
\int\limits_{0}^{x}V(t)|\theta^{2}(t,\lambda)|\,dt\right)(1+o(1)) 
\end{equation}
as $y\to\infty$.

\rm Both theorems will follow from a certain general criterion. 
Suppose that $H_{U}$ is an operator for which all solutions of the equation 
\begin{equation}
 -y'' + U(x)y = \lambda y
\end{equation}
are bounded for almost every $\lambda \in S,$ where $S$ is a certain set 
of positive
 Lebesgue measure. It is known (see, e.g.\ \cite{Sim1}, \cite{Sto}) that in 
this case the set $S$ belongs to an essential support of the absolutely 
continuous part of the spectral measure. 
 Let us pick a family of solutions $\theta(x,\lambda),$ $\lambda \in S,$
 of the equation (3), such that $\theta (x,\lambda)$ are uniformly bounded
 over $S$ and $\theta(x,\lambda),$ $\overline{\theta }(x,\lambda)$ are 
linearly independent for every $\lambda \in S.$ It is easy to see that 
we can always find such family. We have \\

\noindent  \bf Theorem 1.3. \it  Suppose that the potential $V(x)$ is such 
that there exist $\epsilon>0$ and $p\leq 2$ 
so that $(1+x)^{\epsilon}V(x)\in L^{p}.$
 Assume that there exist measurable functions $\theta(x,\lambda)$ 
satisfying the above conditions, such that the operator
\begin{equation}
(Kf)(\lambda) = \chi(S)\int\limits_{0}^{\infty}\theta(x,\lambda)^{2} \exp 
\left(\frac{i}{\Im(\theta\overline{\theta}')}\int\limits_{0}^{x}V(t)
|\theta(t,\lambda)|^{2}\,dt\right)f(x)\,dx 
\end{equation}
satisfies an $L_{2}(R^{+},dx)-L^{2}(S,d\lambda)$ bound on functions $f$ of 
compact support. Then the absolutely continuous  spectrum of $H_{U}$ 
supported on the set $S$ is 
preserved under perturbation by $V$, that is, the set $S$ belongs to an 
essential support of the absolutely continuous part of the spectral 
measure of operator $H_{U+V}.$
Moreover, for almost every $\lambda \in S,$ there exist solutions 
$\psi_{\lambda}(x),$ $\overline{\psi}_{\lambda}(x)$ of the equation (1)
with the asymptotic behavior (2).  \\

\noindent \it Remarks. \rm 1. The assumption 
of the boundedness of all solutions at almost all energies corresponding 
to the essential support of absolutely continuous spectrum is 
rather natural. Almost all known examples of one-dimensional Schr\"odinger
 operators with the absolutely continuous spectrum satisfy 
this assumption. Only recently there appeared rigorous counterexamples to 
the conjecture that this is true in general (see \cite{Mol}), but the 
corresponding potentials are of rather special form
and in particular are not bounded from below.

\rm 2. All three main theorems that we prove have 
natural analogues for the whole axis problems. We will not focus on 
this aspect; all proofs may be generalized to the whole axis case 
in a straightforward manner following \cite{Kis1}. 

\rm 3. 
We have not been able to treat potentials that are assumed merely to 
belong to $L^{p}$ for some $1<p \leq 2$. \\

\rm The main new technique we develop in this paper involves norm 
bounds and almost everywhere convergence results for a class
of multilinear integral operators, which may be of interest 
in its own right.
The plan of the paper is as follows. In Section 2 we 
discuss the basic scheme of asymptotic integration. In Section 3 
we formulate first key results on the estimates for maximal functions 
of certain integral operators. In
 Section 4 we establish norm estimates for multilinear integral
 operators, and in Section 5 prove corresponding a.e.\ convergence results.
 Section 6 contains the conclusion of the proof of all main results.
 In the Appendix we discuss generalization of our results to the
case of potentials which may have strong local singularities.

Independently,  results similar to Theorem 1.1
 were obtained by Remling
 \cite{Remac} by a very different method, based in part on the 
ideas from \cite{Kis1} and \cite{Mol1}. Some of the results we prove 
here (along with some of the results of \cite{Remac}) 
were announced in \cite{CKR}. The announcement \cite{CKR} also 
contains a list of open problems  which we find most 
interesting.

\vspace{0.2cm}

\begin{center}
 \large \bf 2. Asymptotic integration and bounded eigenfunctions.
\end{center}

\vspace{0.2cm}

To prove the stability of the absolutely continuous spectrum of the operator
 $H_{U},$ we will use the following Lemma: \\

\noindent \bf Lemma 2.1. \it Let the potential $W(x)$ be locally integrable 
and satisfy
\[ {\rm sup}_{x} \int\limits_{|x-t|\leq 1} |W_{-}(t)|\,dt < \infty \]
($W_{-}$ denotes the negative part of the potential $W$).
Suppose that for every energy $\lambda$ from the set $E$ of positive 
Lebesgue
 measure, all solutions of the equation
\[ -\frac{d^{2}}{dx^{2}}y(x) + W(x)y(x)=\lambda y(x) \]
are bounded. Then the set $E$ belongs to an essential support of the 
absolutely
 continuous part of the spectral measure of the operator $H_{W}.$ \\

\rm The proof of this Lemma may be found in \cite{Sto} (see also 
\cite{Sim3}). We notice that in different formulations, the fact 
that bounded solutions imply absolutely continuous spectrum was known
 for a long time (see, e.g., \cite{Car}). Lemma 2.1 
is the most convenient statement for our purpose. 

The general plan of our proof is similar to \cite{Kis1} and
 may be described as follows. By assumption, 
we know that for every $\lambda \in S,$ all solutions of the generalized 
eigenfunction equation (3) for the unperturbed operator are bounded.
Examples with imbedded eigenvalues \cite{Nab}, \cite{Sim1}, \cite{KRS}
show that if $V(x)$ is not short-range, we cannot hope in general that for
 every $\lambda \in S$ we still have only bounded solutions for a perturbed
 equation; there may exist rather rich, dense in $S$ set for which we will 
have decaying ($L^{2}$) and therefore 
also growing solutions. Our goal will be to show that nevertheless 
for a.e. $\lambda \in S,$ we still have only bounded solutions. This will 
ensure that the absolutely continuous spectrum is preserved, although 
embedded singular spectrum may occur.

Thus, our goal is to study the solutions of the equation 
\[ -y''+(U(x)+V(x))y=\lambda y. \]
We rewrite this
equation as a system
\[ y_{1}' = \left( \begin{array}{cc} 0 & 1 \\ U+V-\lambda &
 0 \end{array} \right)y_{1}, \]
where $y_{1}$ is now a vector $\left(
 \begin{array}{c} y \\ y' \end{array} \right).$
Let us apply a variation of the parameters transformation with solutions 
of the unperturbed equation
\[ y_{1} = \left( \begin{array}{cc} \theta (x, \lambda) & \overline{\theta}
 (x, \lambda)
\\ \theta '(x, \lambda) & \overline{\theta}' (x, \lambda) \end{array} 
\right) y_{2}, \]
to bring the equation to a more symmetric form
\begin{equation} y_{2}' = \frac{i}{2\Im (\theta \overline{\theta}')}\left( 
\begin{array}{cc} V(x)|\theta (x, \lambda)|^{2} & V(x) \overline{\theta}
 (x, \lambda)^{2} \\ -V(x) \theta (x, \lambda)^{2} & -V(x) |\theta (x, 
\lambda)|^{2} \end{array} \right)y_{2}. 
\end{equation}
Notice that $2i\Im (\theta, \overline{\theta'})$ is a Wronskian of two 
solutions $\theta$ and $\overline{\theta}$ and hence is independent of $x.$
Let us introduce the notation
\[ p(x,\lambda) = \frac{1}{2\Im (\theta \overline{\theta'})} 
\int\limits_{0}^{x}V(y) |\theta(y,\lambda)|^{2}\,dy. \]
It will be convenient to apply to (5) an additional transformation
\[ y_{2} = \left( \begin{array}{cc} \exp (ip(x,\lambda))  & 0
 \\ 0 & \exp (-ip(x,\lambda)) \end{array} \right)y_{3} . \]
We arrive at the following equation for $y_{3}:$
\begin{equation}
y_{3}' = \frac{i}{2\Im (\theta \overline{\theta}')}\left( 
\begin{array}{cc} 0 & V(x) \overline{\theta}
 (x, \lambda)^{2}\exp (2ip(x,\lambda)) \\ -V(x) \theta (x, \lambda)^{2}
\exp (-2ip(x,\lambda)) & 0 \end{array} \right)y_{3}. 
\end{equation}
We  follow the idea of ``$(I+Q)$" asymptotic integration 
originating in Harris-Lutz \cite{HaLu}: to find some invertible
 transformation of equation (6) which would make off-diagonal 
terms absolutely integrable and then  apply 
Levinson's theorem (see, e.g. \cite{CoLe}) to find the asymptotic 
behavior of solutions of the resulting equation. If we succeed, 
we can go back and find the asymptotic behavior of the solutions
 of the original equation.

Let
\begin{equation}
 y_{3} = (1-|q|^{2})^{-\frac{1}{2}}\left( I+ Q \right) y_{4}, 
\end{equation}
where $I$ is an identity matrix, 
\[ Q = \left( \begin{array}{cc} 0 & q \\ \overline{q} & 0 \end{array} 
\right), \]
and $q(x,\lambda)$ is some function to be defined. 
A computation gives for $y_{4}$
\begin{eqnarray}
y_{4}' & = & 
(1-|q|^{2})^{-1}  \frac{i}{2\Im (\theta \overline{\theta}')}
\Bigg[
\nonumber \\
&& 
\left( \begin{array}{cc} V \Re (\theta^{2}q \exp(-2ip)) & V
\overline{\theta}^{2}\exp(2ip) + V\theta^{2}q^{2}\exp(-2ip) \\
- V\theta^{2}\exp(-2ip) - V\overline{\theta}^{2}\overline{q}^{2}\exp(2ip) 
& -V \Re (\theta^{2}q \exp(-2ip)) \end{array} \right) 
\ y_4
\nonumber
\\
&& +\qquad  
\left(  \begin{array}{cc} \frac{1}{2}
(q\overline{q}'
-q'\overline{q}) &
-q' \\ -\overline{q}' & \frac{1}{2}(q'\overline{q}-q\overline{q}') 
\end{array}
 \right) \  y_{4} \Bigg].
\end{eqnarray}
We summarize the main result of this section in \\

\noindent \bf Theorem 2.2. \it Suppose that for some given $\lambda$ there 
exists a function $q(x,\lambda) \in C^{1},$ such that
$q(x,\lambda) \to 0$ as $x\to +\infty$, such that 
\begin{eqnarray}
q'(x,\lambda) & + & \frac{i}{2\Im(\theta \overline{\theta})}V(x)\overline
{\theta}^{2}(x,\lambda)\exp (2ip(x,\lambda)) 
\nonumber \\
&& \qquad + \frac{i}{2\Im(\theta 
\overline{\theta})}V(x) \theta(x,\lambda)^{2}\exp(-2ip(x,\lambda))
q^{2}(x,\lambda) \qquad \in L^{1}.
\end{eqnarray}
Then all solutions of the generalized eigenfunction equation (1) are
 bounded. Moreover, there are two solutions $\psi(x, \lambda),$ 
$\overline{\psi}(x,\lambda)$ with the asymptotic behavior 
\begin{eqnarray}
\psi(x,\lambda) 
& = & \theta(x,\lambda) \exp (i p(x,\lambda)) 
\nonumber \\
&& \cdot \ \exp \left(
\frac{i}{2\Im(\theta \overline{\theta})} \int\limits_{0}^{x}
(1-|q|^{2})^{-1}V\Re(\theta^{2}q\exp(-2ip))\,dt \right)
\left(1+o(1)\right).
\end{eqnarray}

\noindent \it Remark. \rm We note that in all applications that we will 
 have the last 
cumbersome term in the product giving the asymptotic behavior will turn 
out to be integrable. \\

\begin{proof} \rm The proof follows immediately from Levinson's 
theorem (see, e.g. \cite{CoLe}) and equation (8). We may consider
 this system of equations only for $x$ large enough, so that 
$|q(t,\lambda)|<1$ for all $t>x$ and the transformation (7) is non-singular.
 By the assumption of the theorem, the off-diagonal terms
 are absolutely integrable. The diagonal terms are purely imaginary and 
hence
 Levinson's theorem is applicable. Asymptotic behavior of the solutions (10) 
follows directly from the explicit solution of the equation (8) with
 diagonal terms omitted and application of 
transformations we applied to the original system of equations. 
\end{proof} \rm

To complete the proof of Theorem 1.3, we need to construct the function
 $q(x,\lambda)$ verifying (9) and conditions given in Theorem 2.2 for  
almost every
$\lambda \in S$. The main problem is that if we try to solve the equation 
\begin{eqnarray}
q'(x,\lambda) & + &
\frac{i}{2\Im(\theta \overline{\theta})}V(x)
\overline{\theta}^{2}(x,\lambda)\exp (2ip(x,\lambda)) 
\nonumber \\ 
&& \qquad +\ \frac{i}{2\Im(\theta \overline{\theta})}V(x) \theta^{2}(x,\lambda)
\exp(-2ip(x,\lambda))q^{2}(x,\lambda) =0
\end{eqnarray}
by iteration, we obtain  expressions involving multilinear integral
 operators of certain type. We need to show that these expressions 
converge for a.e.\ $\lambda \in S$ in order to ensure $q(x,\lambda) 
\stackrel{x \rightarrow \infty}{\longrightarrow} 0$ for a.e.\ 
$\lambda,$ and to make sure that (9) is satisfied after some
number of iterations.
The first approximation to the solution would be 
\[ q^{(0)}(x,\lambda) = \frac{i}{2\Im(\theta \overline{\theta})}
\int\limits_{x}^{\infty}V(t)\overline{\theta}^{2}(t,\lambda)
\exp (2ip(t,\lambda))\,dt. \]
Again, we have to justify this formula by proving that the 
conditional integral is well-defined for almost every $\lambda.$ This is
 relatively simple  and has been already done in \cite{Kis1}.
 In the next section, we formulate the main result from \cite{Kis1} that we
 will use and make a few comments.


\vspace{0.2cm}

\begin{center}
 \large \bf 3. Almost everywhere convergence for integral operators.
\end{center}

\vspace{0.2cm}

Let an operator $K$ be defined on all measurable bounded functions $f$
of compact support by
\begin{equation}
(Kf)(\lambda) = \int\limits_{0}^{\infty} k(\lambda,x) f(x)\, dx,
\end{equation}
where $k(\lambda,x)$ is a measurable and bounded function on 
$I \times R^{+}.$
To study the a.e.\ convergence of the integral defining $Kf(\lambda)$
on functions from $L^{p},$ we study the corresponding maximal function.
Denote by $M_{K}f(\lambda)$ the maximal function
\begin{equation}
M_{K}f(\lambda) = {\rm sup}_{N} \left| 
\int\limits_{0}^{N}k(\lambda,x)f(x)\,dx
 \right|.
\end{equation}
The following is a mild generalization of a result proved in
\cite{Kis1}. \\ 

\noindent \bf Theorem 3.1. \it 
Let $p,q$ be exponents satisfying $1\le p<q \le \infty$. 
Suppose that $K$ is a bounded linear operator from $L^p(R)$
to $L^q(R)$. 
Then the maximal function $M_K$ also maps $L^p$ to $L^q$ boundedly,
that is,
\begin{equation}
 \|M_{K}f\|_{q} \leq C_{q} \|f\|_{p} \;\: for \;\: every \;\: f \in L^{p}
\end{equation}
As a consequence, the integral 
\[ \int^{N}_{0} k(\lambda,x)f(x)\,dx \]
converges
as $N \rightarrow \infty$ for almost every value of 
$\lambda$, for every $f \in L^{p}.$  \\

\rm For the proof of a slightly less general result
we refer to \cite{Kis1}; we will also 
sketch in the appendix a proof of a very similar result 
(Lemma A.3). Theorem 3.1 may be obtained from that
proof by a simple modification.

The following variant may also be proved by the same method. 
Denote by $\|K\|_{p,q}$ the norm of $K$, as an operator
from $L^p$ to $L^q$. \\

\noindent \bf Theorem 3.2. \it 
Let $(X,\mu)$, $(Y,\nu)$ be $\sigma$-finite measure spaces.
Suppose that $1\le p<q\le\infty$.
Then for any bounded linear operator $K$ from $L^p(X)$ to $L^q(Y)$, 
and for any sequence of measurable
sets $\{E_n\subset X: n\in {Z}\}$ such that $E_n \subset
E_{n+1}$ for every $n$,
the maximal function
\[
M_K f(y) = \sup_n |K(f\cdot\chi_{E_n})(y)|
\]
is bounded from $L^p(X)$ to $L^q(Y)$.
Moreover
\[
\|M_K\|_{p,q} \le A\|K\|_{p,q},
\]
where $A<\infty$ depends only
on $p,q$. \\

\rm
See \cite{Kis2} for other results along these lines.

\vspace{0.2cm}

\begin{center}
\large \bf 4. Norm estimates for multilinear transforms.
\end{center}

\vspace{0.2cm}

In this section, we study the questions related to the norm estimates  
for certain multilinear transforms. The results of this section and 
the next 
 will enable us to fulfill the plan sketched in the end of Section 2 and
 find the function $q(x,\lambda)$ with the needed properties
 for a.e.\ $\lambda$ by iteration of (11). 

Suppose that the functions $k_{i}(\lambda, x),$ $i=1,n \dots $ are 
defined on $I \times R^{+},$ where $I$ is some measurable set in $R.$ 
We assume that the operators 
\[ (K_{i}f)(\lambda) = \int\limits_{R^{+}} k_{i}(\lambda, x)f(x)\,dx \]
satisfy the bounds 
\begin{equation}
 \|K_{i}f\|_{L^{q}(I, d\lambda)} \leq C_{i}\|f\|_{L^{p}(R^{+}, dx)} 
\end{equation}
on functions of compact support
for some $2>p \geq 1$ and $q>p.$

Let $n \geq 2$. Let $A$ be 
any set of ordered  pairs $\alpha = (i_{\alpha}, i'_{\alpha}),$
with $1 \leq i_{\alpha}, i'_{\alpha} \leq n.$ Let $|A|$ denote 
the cardinality of $A.$ By $\chi_{E}(x)$ we denote 
a characteristic function which is equal to one  when $x \in E$ and  
is zero otherwise.

Consider the multilinear operator $T_{n}$ given by 
\begin{equation}
T_{n}(f_{1}, \dots f_{n})(\lambda) 
= \int\limits_{R^{n}} \prod\limits_{j=1}^{n}
f_{j}(x_{j})k_{j}(x_{j},\lambda) \prod\limits_{\alpha \in A}
\chi_{R^{+}}(x_{i_{\alpha}}-x_{i'_{\alpha}})\,dx, 
\end{equation}
$x=(x_{1}, \dots x_{n}).$ Notice that if there were no ``diagonal" 
characteristic functions, 
the expression (16) would decompose into a product of one-dimensional 
integrals, and the analysis would become trivial. \\

\noindent \it Remark. \rm  We do not rule out the possibility
that some of the characteristic functions in (16) are contradictory 
and the whole expression is zero. \\

Our goal in this section is to prove the following 
property: \\

\noindent \bf Theorem 4.1. \it Suppose that the multilinear operator
 $T_{n}$ is given by (16) with kernels $k_{j}(\lambda,x_{j})$ 
satisfying (15).  
Then for any functions $f_{i} \in L^{p}(R^{+},dx),$ $i=1, \dots n,$
 such that the integral (16)
 converges absolutely for a.e.\ $\lambda,$ we have 
\[ \|T_{n}(f_{1},...f_{n})\|_{s_{n}} \leq C_{n}\prod_{i=1}^{n}
\|f_{i}\|_{p}, \]
where $s_{n}^{-1} = nq^{-1}.$ The constant $C_{n}$ depends only on $n$ and
 constants in the norm bounds (15) for operators $K_{i}.$  \\

\noindent \it Remarks. \rm The conclusion of Theorem 4.1 holds 
in particular when $s_{n}<1.$ Our proof will yield a
more general inequality,
in which $f_j\in L^{p_j}$, the exponents $p_j$ vary freely in 
a certain range, and $s_n^{-1} = \sum_j q_j^{-1}$ where
$q_j^{-1} = 1 - p_j^{-1}$; the case where 
all exponents $p_j$ are equal suffices for our applications,
and we restrict attention to it in order to simplify 
computations slightly.
\\

\noindent \rm  By assumption, the value of $T_{n}(f_{1}, \dots
f_{n})(\lambda)=
g(\lambda)$ is well-defined for a.e.\ $\lambda$ by the absolutely convergent 
integral.
Our strategy will be to divide the domain of integration into 
disjoint pieces
 and represent the function $g$ as a sum of terms coming from 
integration over
 these disjoint pieces, formally:
 \[ g(\lambda)=\sum_{i=1}^{\infty}g_{i}(\lambda).\]
Because of the absolute convergence, we have that the sum 
$\sum_{i=1}^{n}g_{i}
(\lambda)$ converges to $g(\lambda)$ for a.e.\ $\lambda$ as $n \rightarrow 
\infty.$ We show, choosing the functions $g_{i}$ in a convenient way, that
 the sum also converges absolutely in the 
appropriate space $L^{s_{n}},$ thus proving Theorem 4.1.

In the proof of Theorem 4.1, we will need a certain representation of the 
function 
\[ f(x_{1})f(x_{2})\chi_{R^{+}}(x_{2} - x_{1}) \]
 as a sum of products of two functions depending only on $x_{1}$ 
and $x_{2}$ 
respectively.
Let us first introduce a decomposition of $R^{+}$ associated with 
the function
 $f.$ Normalize the function $f$ so that  $\|f\|_{p}^{p}=1.$  
By $\chi_{E}$ we
 will denote the characteristic function of the set $E.$ 
Let $E(1,1)$ and
 $E(1,2)$ be 
disjoint intervals such that 
\[ \|f(x)\chi_{E(1,1)} \|^{p}_{p} = 
\|f(x)\chi_{E(1,2)} \|^{p}_{p} =
 2^{-1}, \]
$E(1,1) \cup E(1,2)=R^{+}$ and $E(1,1)$ lies 
entirely to the 
right of $E(1,2)$ (i.e. for any $x \in E(1,1),$ 
$y \in E(1,2)$
 we have $x \leq y$). We note that $E(1,2)$ is 
half-infinite and assume 
$E(1,1)$ contains its 
right end for the above decomposition to hold. We also remark that the 
decomposition is not necessarily unique ($f$ might vanish on some set so 
that this decomposition will be non-unique), and we just take some 
decomposition. In future we will omit such
 inessential details. We continue to decompose each of the intervals
 $E(1,l)$ in a similar manner, obtaining on the $m^{{\rm th}}$ step
  $2^{m}$
  intervals $\{ E(m,l)\}_{l=1}^{2^{m}},$ such that 
$\cup_{l=1}^{2^{m}}E(m,l)=R^{+},$ $\|f(x)\chi (E(m,l))
\|_{p}^{p}=2^{-m}$ for $j=1,...2^{m},$ the intervals are disjoint
 and $E(m,l)$ lies entirely to the left from $E(m,i)$ if
 $l < i.$ In notation $E(m,l),$ we refer to $m$ as ``generation"
 of this interval and to $l$ as ``index". Of 
importance, in particular, will be the following  evident property of
 intervals $\{E(m,l): 1 \le m \le \infty, 1 \le l \le 2^{m} \}:$
 any two intervals are either disjoint or one is contained in another. 

We  proceed to decompose the ``diagonal" characteristic functions
 in a convenient way. 

\begin{figure}
\includegraphics{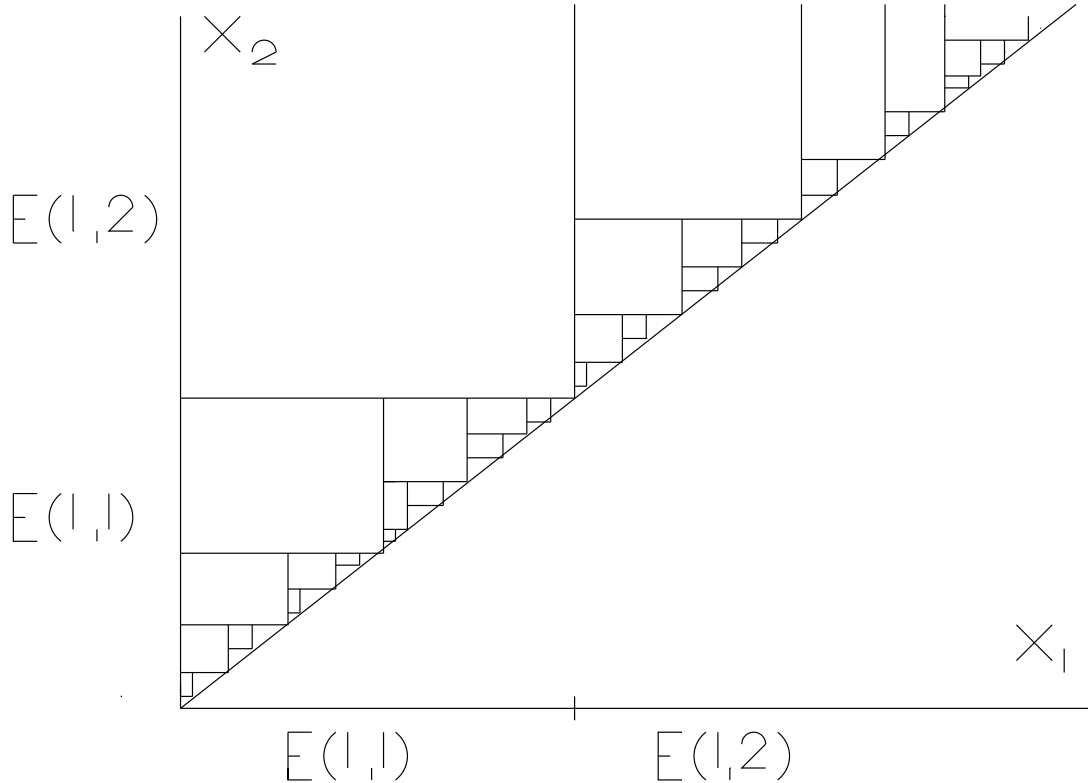}
\caption{Decomposition of $\chi_{R^{+}}(x_{2}-x_{1})$}\label{Fi:A}
\end{figure}

\noindent \bf Lemma 4.2. \it The following identity holds:
\begin{equation}
\chi_{R^{+}}(x_{2} - x_{1}) f(x_{1})f(x_{2})= 
\left( \sum\limits_{m=1}^{\infty}
\sum\limits_{ 
\stackrel{\scriptstyle l=1} {\scriptstyle l \,{\rm odd}}}
^{2^{m}} 
\chi_{E(m,l)}(x_{1}) 
\chi_{E(m,l+1)}(x_{2}) \right)f(x_{1})f(x_{2}).  
\end{equation}

\begin{proof} \rm Figure 1 illustrates the 
decomposition that we perform. Let us denote by $H_{12}$ the set
\[ H_{12} = \{ x \in R^{2}, x=(x_{1},x_{2}) | x_{1} < x_{2} \}, \]
and by ${\rm supp} f$ the closure of the set of the points $x$ such
that for every interval $I,$ such that $x \in I,$ $|f(x)|$ is positive
 on the set of positive Lebesgue measure in $I.$
The claim will follow if we show that 
\[ H_{12} \cap ({\rm supp}_{x_{1}}f \times {\rm supp}_{x_{2}}f)=
 \bigcup_{m=1}^{\infty} 
\bigcup^{2^m}_{
\stackrel{\scriptstyle l=1} {\scriptstyle l \,{\rm odd}}}
(E(m,l) 
\times E(m,l+1) 
) \]
and the sets under the union on the right hand side are disjoint. 
The latter fact is easy to see: if $E(m,l) \subset E(s,i),$
 $l$ odd, $s\neq m,$ then necessarily $m>s$ and $E(m,l+1)$ also
 belongs to $E(s,i),$ not $E(s,i+1).$ On the other hand, we
 show that for every $y_{2},y_{1} \in {\rm supp}
 f,$ $y_{1} < y_{2},$ there exist two sets $E(m,l),$ $E(m,l+1),$ 
with $l$ odd, such that $y_{1} \in E(m,l)$ and $y_{2}
 \in E(m,l+1).$
 Let $\|f \chi_{(y_{1}, y_{2})}\|_{p}^{p}=a>0.$ Here we assume that $f$ is
 normalized and use the 
condition that $y_{1},$ $y_{2}$ lie in ${\rm supp}f$ to infer 
that $a>0.$ Choose
 $s$ so that $2^{-s} \geq a \geq 2^{-s-1}.$ If $y_{1},$ $y_{2}$ lie in one 
set of generation $s,$ $E(s,l),$ then necessarily $y_{1} 
\in E(s+1,2l-1)$ and $y_{2} \in E(s+
1,2l).$ If $y_{1}$ and $y_{2}$ lie in different sets of 
 generation $m$, $E(s,l)$ and $E(s,l+1),$ then either 
$l$ is odd or $y_{1} \in E(s-1,l/2^{r})$ and 
$y_{2} \in E(s-1,l/2^{r}+1),$ where $r$ is such that $l/2^{r}$ is 
odd. 
\end{proof} \rm

\noindent \it Remark. \rm In particular, if ${ \rm supp} f = R^{+},$
 we get a representation of diagonal characteristic function $\chi_
{R^{+}}(x_{2} - x_{1})$ as a sum of products of characteristic 
functions of some intervals in $x_{1}$ and $x_{2}$ variables.

We now begin proof of Theorem 4.1. 

\begin{proof} \rm  Since $T_n$ is multilinear,
we may assume without loss of generality throughout the proof that 
$\|f_{i}\|_{p}^p=\frac{1}{n}$ for all $i=1,...n.$ Let 
\[ f(x) = \left( \sum\limits_{i=1}^{n}|f_{j}(x)|^{p} \right)^{\frac{1}{p}}. 
\]
Consider the family of the intervals $\{E(m,l)\}$ associated with the 
function $f.$ An important property of this family is that 
\begin{equation}
 \|f_{i}(x)\chi_{E(m,l)}\|^{p}_{p} \leq 2^{-m} 
\end{equation}
for all $i,l,m$.
Write 
\[ A=\{ \alpha_{1}, \dots \alpha_{|A|} \}. \]
We begin by substituting the result of Lemma 4.2 into formula (16):
\begin{equation}
 T_{n}(\overline{f})= 
 \sum\limits_{m_{1}}^{\infty} \cdots
\sum\limits_{m_{|A|}=1}^{\infty}
 \sum\limits_{l_{1}}' \cdots
\sum\limits_{l_{|A|}}'\,\int\limits_{R^{n}}
dx \prod\limits_{j=1}^{n} k_{j}(\lambda, x_{j})f_{j}(x_{j}) 
\prod\limits_{t=1}^{|A|} \chi_{E(m_{t},l_{t})}(x_{i_{t}})
\chi_{E(m_{t},l_{t}+1)}(x_{i'_{t}}),  
\end{equation}
where $\sum\limits_{l_{t}}'$ means the sum over all odd integers $l_{t}
\in [1,2^{m_{t}}],$ and $i_{t}=i_{\alpha_{t}},$ $i_{t}'=i_{\alpha_{t}}'.$
Thus
\begin{equation}
| T_n(\overline{f})(\lambda)  | \leq 
\sum\limits_{m_{1}}^{\infty} \cdots \sum\limits_{m_{|A|}=1}^{\infty}
F^{\overline{m}}
(\overline{f})(\lambda), 
\end{equation}
where $\overline{m}= (m_{1},m \dots _{|A|})$ and 
\begin{equation}
 F^{\overline{m}}(\overline{f})(\lambda)=
\sum\limits_{l_{1}}' \cdots \sum\limits_{l_{|A|}}' \, \int\limits_{R^{n}}
dx \prod\limits_{j=1}^{n} |K_{j}(f_{j}\chi_{G(j,\overline{l})})(\lambda)|, 
\end{equation} 
where $\overline{l}=(l_{1},...l_{|A|})$ (all variables
 $l_{t}$ take only odd
values), the set $G(j,\overline{l})$ 
depends on $\overline{m}$ (and on $A$), and 
\[ G(j,\overline{l}) = \left[ \bigcap_{t: j=i_{t}}E(m_{t},l_{t}) \right] 
\bigcap
\left[ \bigcap_{t:j=i'_{t}}E(m_{t},l_{t}+1) \right].  \]

We aim to prove that
\begin{equation}
\|F^{\overline{m}}(\overline{g})\|_{s_{n}} \leq C_{n}2^{-\gamma_{n}
|\overline{m}|}\prod\limits_{j=1}^{n}\|g_{j}\|_{p}^{(1-\beta_{n})}
\end{equation}
for any $g_{1},...g_{n}$ which satisfy (18)
(that is, $\| g_i \chi_{E(m,l)}\|_p^p \le 2^{-m}$
for all $i,m,l$).
 Here $\gamma_{n}$ is some positive constant which 
 depends only on $n,$ $C_{n}$ depends on $n$ and constants in norm bounds
(15) for the operators $K_{j},$ and $\beta_{n}$ satisfies 
\[ 1 \geq 1-\beta_{n} > {\rm max}(\frac{p}{q}, \frac{p}{2}) . \]
Invoking the decomposition (20) and summing over $\overline{m}$,
Theorem 4.1 follows directly from (22).

We may proceed to prove (22) by induction on $|A|.$ The case 
$|A|=0$ is immediate from the hypothesis on $K_{j}$ by H\"older's inequality
(in this case $\overline{m}$ is a $0$-tuple). 
It will be convenient to consider the graph 
$\Upsilon$ with vertices $\{1,...n\}$ and edges $(i_{t},i'_{t})$ joining
$i_{t}$ to $i_{t}'$ for any $t$ (and no other edges). 
To each edge we associate the generation $m_{t}$ which 
corresponds to a generation in the decomposition (17)
of $\chi_{R^{+}}(x_{i_{t}}-x_{i_{t}'})$ that we fixed in the sum (21).
It suffices to 
treat the case where $\Upsilon$ is connected; the general case then 
follows by H\"older's inequality. 

Fix $\overline{m}.$ Relabel the indices so that $m_{1} \leq ... \leq 
m_{|A|}.$ For simplicity of notation, we also relabel pairs $(i_{t},i_{t}')$
so that $m_{t}$ still denotes the generation in the decomposition of 
$\chi_{R^{+}}(x_{i_{t}}-x_{i_{t}'}).$
Let $N$ be the largest index for which $m_{N}=m_{1}.$ 

For many values of $\overline{l},$ the set $G(j,\overline{l})$ is 
empty for some $j$.
Such terms contribute zero to the sum (21); 
this observation underlies the estimate (22) 
for $F^{\overline{m}}(\overline{f})$. 
To take this into account,
drop from the sum (21) all terms for which there exists $j$ such that 
$G(j,\overline{l})=\emptyset;$ such terms contribute $0$ 
to $F^{\overline{m}}(\overline{g})$. We say that an index 
$l$ \it remains \rm if the corresponding term has not been dropped.

We have the following \\

\noindent \bf Lemma 4.3. \it For any $1 \leq j \leq n,$ either
\begin{equation}
G(j,\overline{l}) \subset E(m_{1},l_{1}) \,\,\,\forall \overline{l}
 \,\,\,{\rm remaining}
\end{equation}
or 
\begin{equation}
G(j,\overline{l}) \subset E(m_{1},l_{1}+1) \,\,\,\forall \overline{l}
 \,\,\,{\rm remaining}.
\end{equation}
Let $B_{1}=\{j:$ (23) holds $\}$, $B_{2} = \{ j:$ (24) holds $\}$.
Then for any $2 \leq t \leq N$ (that is, if $m_{t}=m_{1}$), 
\begin{equation}
l_{t}=l_{1}\,\,\,\forall \overline{l} \,\,\, {\rm remaining}.
\end{equation}
Finally for each $t>N$ (so $m_{t}>m_{1}$), either $i_{t},$ $i'_{t}$ are 
both in $B_{1}$ for all remaining $l$ 
or they are both in $B_{2}$ for all remaining $l.$
We say that $t \in A_{1},$ $t \in A_{2}$ respectively.  \\

\begin{proof} \rm For any $m \geq 1,$ $l$ odd, set 
\[ \tilde{E}(m,l) = E(m,l) \cup E(m,l+1). \]
First we prove that 
if $\overline{l}$ 
remains, then 
\begin{equation} 
\tilde{E}(m_{t},l_{t})  \subset \tilde{E}(m_{1},l_{1})  
\end{equation}
for all $t.$ 
Notice that both sets in (26) also belong to the 
family $E(m,l)$ (they are $E(m_{t}-1, \frac{l_{t}-1}{2}),$
$E(m_{1}-1, \frac{l_{1}-1}{2})$ respectively; we may 
assume $E(0,0)=R^{+}$). 
Therefore, to prove (26) it is sufficient to show 
that the two sets in (26) intersect, since in this case
one is contained in another by the martingale-type property. 
 
Recall that $m_{1}$ is the generation which is 
fixed  in  decomposition of the  
characteristic function $\chi_{R^{+}}(x_{i_{1}}-x_{i_{1}'})$
in the sum (21) for $F^{\overline{m}}.$
Pick any other $m_{t}$ which is fixed
 in the decomposition of the  characteristic
function $\chi_{R^{+}}(x_{i_{t}}-x_{i_{t}'}).$
Since the graph $\Upsilon$ is connected, we can find 
a path in $\Upsilon$ which connects either $i_{t}$ or $i_{t}'$ with either
$i_{1}$ or $i_{1}',$ and does not contain the edges $(i_{1},i_{1}'),$
$(i_{t},i_{t}').$
 Suppose that this path goes from $i_{1}$ to $i_{t}$ and
passes successively through the edges with 
the corresponding generations $m_{t_{1}},...m_{t_{r}}.$ 
This path does not depend on $\overline{l}$.

For $G(j,\overline{l})$
to be non-zero for all $j,$ we must have 
\[
 \tilde{E}(m_{1},l_{1}) \cap \tilde{E}(m_{t_{1}},l_{t_{1}}) \neq
\emptyset, \,\,\,\, \tilde{E}(m_{t},l_{t}) \cap \tilde{E}(m_{t_{r}},
l_{t_{r}}) \neq \emptyset, \,\,\,\,\,{\rm and} \]
\begin{equation} 
\tilde{E}(m_{t_{i}},l_{t_{i}}) \cap \tilde{E}(m_{t_{i+1}},l_{t_{i+1}})
\neq \emptyset \,\,\,\,{\rm for\,\,\, all\,\,\,} i=1,...r-1. 
\end{equation}
Hence  by our assumption that $m_{1} \leq m_{t}$ for all $t$ we see that 
$\tilde{E}(m_{t_{1}},l_{t_{1}}) \subset \tilde{E}(m_{1},l_{1}).$
But then by (27) also 
$\tilde{E}(m_{1},l_{1}) \cap \tilde{E}(m_{t_{2}},l_{t_{2}}) \neq
\emptyset,$ hence  $\tilde{E}(m_{t_{2}},l_{t_{2}}) 
\subset \tilde{E}(m_{1},l_{1}).$ We continue in the same way 
concluding that 
$\tilde{E}(m_{t},l_{t}) \subset \tilde{E}(m_{1},l_{1})$ and hence (26)
holds. 

The statements (23), (24) and (25)
 of the lemma now follow immediately from the 
martingale-type
property of the sets $E(m,l)$ and the definition of the set 
$G(j,\overline{l}).$
To obtain (25), note that because $l_1$ and $l_t$ are odd,
when $l_t = l_1$
the inclusion $E(m_t,l_t)\subset \tilde E(m_1,l_1)$
forces $E(m_t,l_t) = E(m_1,l_1)$.
To prove the final statement,
 suppose that we know in addition that $m_{t}>m_{1}.$
We can find a path in $\Upsilon$ which goes from $i_{t}$ or $i_{t}'$ 
to a vertex adjacent to an edge $(i_{s},i_{s}')$ 
with the corresponding 
generation $m_{s}$ equal to $m_{1}$ (i.e. with $s \leq N$), and contains only 
edges with the corresponding generations strictly less than $m_{1}.$  
An argument analogous to the above shows that in this case
$\tilde{E}(m_{t},l_{t})$ is contained either in  $E(m_{1},l_{s})$ or in 
$E(m_{1},l_{s}+1)$ for all 
remaining $\overline{l},$ 
depending on whether the vertex to which the path
leads coincides with $i_{s}$ or $i_{s}'$ respectively. By (25)
 the lemma is proven.
\end{proof} \rm 

\noindent \it Remark. \rm It may happen that there exist two (or more)
different paths from $i_{t}$ (or $i_{t}'$),
one of which leads to a vertex
$j_{1}$ where $f_{j_{1}}(x_{j_{1}})$ is multiplied by 
$\chi_{E(m_{1},l_{1})}(x_{j_{1}})$,
while the other leads to a vertex $j_{2}$ where 
$f_{j_{2}}(x_{j_{2}})$ is multiplied by 
$\chi_{E(m_{1},l_{1}+1)}(x_{j_{2}})$. (The simplest way for this
to happen is for $A$ to
include two pairs $(i,j)$ and $(j,i)$.) In this case, $G(j,\overline{l})$
is zero and hence $\overline{l}$ does not remain. \\

By Lemma 4.3
\begin{eqnarray*} 
F^{\overline{m}}(\overline{f})(\lambda)  
& = & \left( \sum\limits_{l_{1}} '
\sum\limits_{l_{N+1}} '\cdots \sum\limits_{l_{|A|}} '
\prod\limits_{j \in B_{1}}|K_{j}(f_{j}\chi_{E(m_{1},l_{1})}\chi_{G(j,
\overline{l})})(\lambda)| 
\right. \\ & & \cdot \ \left.
\prod\limits_{j \in B_{2}}|K_{j}(f_{j}\chi_{E(m_{1},l_{1}+1)}\chi_{G(j,
\overline{l})})(\lambda)| \right)
\\
& \leq &  \sum\limits_{l_{1}} 'F_{A_{1},B_{1}}^{\overline{m}^{(1)},l_{1}}
(\overline{f})(\lambda)F_{A_{2},B_{2}}^{\overline{m}^{(2)},l_{1}}
(\overline{f})(\lambda), 
\end{eqnarray*}
where 
\[ F_{A_{1},B_{1}}^{\overline{m}^{(1)},l_{1}}(\overline{f})(\lambda)=
\sum\limits_{l_{t}:t \in A_{1}} ' \prod\limits_{j \in B_{1}}
|K_{j}(f_{j}\chi_{E(m_{1},l_{1})}\chi_{G(j,l)})(\lambda)| \]
and $\sum\limits_{l_{t}:t \in A_{1}} '$ denotes the sum over all $l_{t}$
such that $1 \leq l_{t} \leq 2^{m_{t}},$ $l_{t}$ is odd, $t \in A_{1},$
 and we write $\overline{m}^{(1)}=
(m_{t})_{t \in A_{1}}.$ Note that 
$F_{A_{1},B_{1}}^{\overline{m}^{(1)},l_{1}}(\overline{f})(\lambda)$
depends only on those $f_{j}$ for which $j \in B_{1}.$ The factor
$F_{A_{2},B_{2}}^{\overline{m}^{(2)},l_{1}}(\overline{f})(\lambda)$
is defined similarly, but $\chi_{E(m_{1},l_{1})}$ is replaced by 
$\chi_{E(m_{1},l_{1}+1)}$. 

We may rewrite for each $j \in B_{1}$ 
\[ G(j,l) = E(m_{1},l_{1}) \bigcap G^{1}(j, (l_{t})_{t \in A_{1}}), \]
where 
\[ G^{1}(j, (l_{t})_{t \in A_{1}}) = \left[ \bigcap\limits_{t \in A_{1},
j=i_{t}} E(m_{t},l_{t}) \right] \bigcap \left[ \bigcap\limits_{t \in A_{1},
j=i_{t}'}E(m_{t},l_{t}+1) \right]. \]
Indeed, by Lemma 4.3 all other sets that enter in the definition 
of $G(j,l)$ belong to $E(m_{1},l_{1}+1)$ and hence are absent
for $\overline{l}$ which remain.
Thus $F_{A_{1},B_{1}}^{\overline{m}^{(1)},l_{1}}(\overline{f})(\lambda)$
and $F_{A_{2},B_{2}}^{\overline{m}^{(2)},l_{1}}(\overline{f})(\lambda)$
are expressions of the same form as the original $F^{\overline{m}}.$
Since $0<|A_{1}|<|A|,$ both 
$F_{A_{1},B_{1}}^{\overline{m}^{(1)},l_{1}}(\overline{f})(\lambda)$ and
$F_{A_{2},B_{2}}^{\overline{m}^{(2)},l_{1}}(\overline{f})(\lambda)$
may be estimated by induction on $|A|.$ Therefore
\begin{equation}
 \|F_{A_{1},B_{1}}^{\overline{m}^{(1)},l_{1}}(\overline{f})(\lambda)\|_
{s_{|B_{1}|}} \leq C2^{-\gamma_{|B_{1}|} |m^{(1)}|} 
\prod\limits_{j \in B_{1}}
\|f_{j}\chi_{E(m_{1},l_{1})}\|_{p}^{1-\beta_{|B_{1}|}}; 
\end{equation}
a similar bound also holds for
$F_{A_{2},B_{2}}^{\overline{m}^{(2)},l_{1}}(\overline{f})(\lambda)$.

Using (28), we are ready to estimate 
$\|F^{\overline{m}}(\overline{f})\|_{s_{n}}.$ We distinguish between two
cases: $s_{n} \le 1$ and $s_{n} \geq 1.$ Suppose first that $s_{n} \le 1$.
Then 
\begin{equation}
 \|F^{\overline{m}}(\overline{g})\|_{s_{n}}^{s_{n}} \leq
\sum\limits_{l_{1}}' 
\|F_{A_{1},B_{1}}^{\overline{m}^{(1)},l_{1}}(\overline{g})(\lambda)\|
_{s_{|B_{1}|}}^{s_{n}}
\|F_{A_{2},B_{2}}^{\overline{m}^{(2)},l_{1}}(\overline{g})(\lambda)\|_
{s_{|B_{2}|}}^{s_{n}}. 
\end{equation}
We used the fact that $\| 
\sum h_{i}(x) \|^{s}_{s} \leq \sum \|h_{i}(x)\|_{s}^{s}$ when $s<1$ and 
H\"older's inequality. Plugging the estimate (28)
 and a similar 
bound for 
$F_{A_{2},B_{2}}^{\overline{m}^{(2)},l_{1}}(\overline{g})(\lambda)$
into (29),  we find
\begin{eqnarray}
 \|F^{\overline{m}}(\overline{g})\|_{s_{n}}^{s_{n}} 
& \leq &
C_{n} 2^{-(\gamma_{|B_{1}|}|m^{(1)}|+\gamma_{|B_{2}|}|m^{(2)}|)s_{n}}
\nonumber
\\
& & \cdot \ \sum\limits_{l_{1}} '\left(
\prod\limits_{j \in B_{1}} \|g_{j}\chi_{E(m_{1},l_{1})}
\|_{p}^{(1-\beta_{|B_{1}|})s_{n}} \prod\limits_{j \in B_{2}}\|g_{j}
\chi_{E(m_{1},l_{1}+1)}\|_{p}^{(1-\beta_{|B_{2}|})s_{n}}\right). 
\end{eqnarray}

Pick  $0 < a_{1}, a_{2} <1$ such that 
\[ a_{1}(1-\beta_{|B_{1}|}) = a_{2} (1-\beta_{|B_{2}|}) = \frac{p}{q}. \]
We can find such $a_{1},$ $a_{2}$ by the induction assumption. 
The sum in (30) may be estimated by H\"older's inequality in the following 
way:
\begin{eqnarray*}
\lefteqn{
\sum\limits_{l_{1}}'
\left(\prod\limits_{j \in B_{1}} \|g_{j}\chi_{E(m_{1},l_{1})}
\|_{p}^{(1-\beta_{|B_{1}|})s_{n}} 
\prod\limits_{j \in B_{2}}\|g_{j}
\chi_{E(m_{1},l_{1}+1)}\|_{p}^{(1-\beta_{|B_{2}|})s_{n}}\right) 
}
\\
& &
\qquad
\le
\prod\limits_{j \in B_{1}}
{\rm max}_{l_{1}}\|g_{j}\chi_{E(m_{1},l_{1})}\|_{p}^
{(1-a_{1})(1-\beta_{|B_{1}|})s_{n}}
\nonumber \\
& & \qquad\qquad \cdot \ \ \prod\limits_{j \in B_{2}}{\rm 
max}_{l_{1}}\|g_{j}\chi_{E(m_{1},l_{1}+1)}\|_{p}^
{(1-a_{2})(1-\beta_{|B_{2}|})s_{n}}\prod\limits_{j=1}^{n}\|g_{j}\|_{p}^
{\frac{p}{q}s_{n}}. 
\end{eqnarray*}
Thus
\begin{eqnarray*}
 \|F^{\overline{m}}(\overline{g})\|_{s_{n}} 
& \leq & C_{n}2^{-{\rm min}(\gamma_
{|B_{1}|}, \gamma_{|B_{2}|})(|m^{(1)}|+|m^{(2)}|)} \\
& & \cdot \ \ \prod\limits_{j=1}^{n}
{\rm sup}_{l}\|g_{j}\chi_{E(m_{1},l)}\|_{p}^{{\rm min}(1-\beta_{|B_{1}|},
1-\beta_{|B_{2}|})-\frac{p}{q}}\prod\limits_{j=1}^{n} \|g_{j}\|_{p}^
{\frac{p}{q}}. 
\end{eqnarray*}
Obviously $\|g_{j}\|_{p} \geq {\rm sup}_{l}\|g_{j}\chi_{E(m_{1},l)}\|_{p},$ 
and $ {\rm sup}_{l}\|g_{j}\chi_{E(m_{1},l)}\|_{p} \leq 2^{-m_{1}}$ by (18).
Pick $\beta_{n}$ so that 
\[ {\rm min}(1-\beta_{|B_{1}|}, 1-\beta_{|B_{2}|}) > 1-\beta_{n} > {\rm 
max}(\frac{p}{q}, \frac{p}{2}), \]
and 
\[ \gamma_{n} = {\rm min}\big(
[{\rm min}(1-\beta_{|B_{1}|},1-\beta_{|B_{2}|})-
(1-\beta_{n})], \gamma_{|B_{1}|}, \gamma_{|B_{2}|}\big). \]
Then 
\begin{equation}
 \|F^{\overline{m}}(\overline{g})\|_{s_{n}} \leq C_{n} 2^{-\gamma_{n}|m|}
\prod\limits_{j=1}^{n} \|g_{j}\|_{p}^{1-\beta_{n}}. 
\end{equation}
There are only finitely many pairs $B_{1},$ $B_{2}$ such that $|B_{1}|+
|B_{2}|=n,$ and hence the constants $\gamma_{n},$ $\beta_{n}$ may be chosen
to be independent of $B_1,B_2$. 

The case $s_{n} \ge 1$ is similar. Using the triangle inequality 
and H\"older's inequality, we get
\begin{eqnarray*} 
\|F^{\overline{m}}(\overline{g})\|_{s_{n}} 
& \le &
\sum\limits_{l_{1}}' 
\|F_{A_{1},B_{1}}^{\overline{m}^{(1)},l_{1}}(\overline{g})(\lambda)\|
_{s_{|B_{1}|}}
\|F_{A_{2},B_{2}}^{\overline{m}^{(2)},l_{1}}(\overline{g})(\lambda)\|_
{s_{|B_{2}|}} 
\\
& \leq &
C_{n}2^{-{\rm min}(\gamma_{|B_{1}|}, \gamma_{|B_{2}|})(|m^{(1)}|
+|m^{(2)}|)}
\\
&& \cdot\ \sum\limits'_{l_{1}} \left(\prod\limits_{j \in B_{1}} 
\|g_{j}\chi_{E(m_{1},l_{1})}
\|_{p}^{(1-\beta_{|B_{1}|})} \prod\limits_{j \in B_{2}}\|g_{j}
\chi_{E(m_{1},l_{1}+1)}\|_{p}^{(1-\beta_{|B_{2}|})}\right). 
\end{eqnarray*}
Provided that 
\begin{equation}
 \frac{1}{p} (|B_{1}|(1-\beta_{|B_{1}|})+|B_{2}|(1-\beta_{|B_{2}|}))>1, 
\end{equation}
we can apply the same argument as in the case $s_{n}<1$ to prove (31).
But (32) holds for all $|B_{1}|,$ $|B_{2}|\geq 1$ since by induction 
hypothesis $1-\beta_{r}> \frac{p}{2}$ for all $r <n.$ This completes 
the proof of (22), and hence of Theorem 4.1. 
\end{proof} \rm


\vspace{0.2cm}

\begin{center}
\large \bf 5. Almost everywhere convergence of multilinear transforms 
\end{center} 

\vspace{0.2cm}

In the proof of the a.e. convergence, an important role will be played by 
the following operators. Let 
$D_{1}(\lambda),...D_{n}(\lambda)$
 be measurable functions of $\lambda$ 
mapping $I$ to $R^{+}\cup\{\infty\}.$ Let us denote by
\[ T^{D_{1}(\lambda),...D_{n}(\lambda)}_{n}(f_{1},...f_{n})
(\lambda) \]
an operator obtained from $T_{n}$ by replacing the
 kernels $k_{i}(\lambda,x_{i}),$ $i=1,...n$ with 
\[
 \tilde{k}_{i}(\lambda, x_{i}) = k_{i}(\lambda, x_{i})
\chi_{R^{+}}(D_{i}(\lambda)-x_{i}), \] 
$i=2,...n.$ Throughout this section, we will assume that the kernels
$k_{i}(\lambda,x)$ are bounded and the
integral operators $K_{i}$ corresponding to these kernels satisfy
$L^{2}-L^{2}$ estimates. The results we prove here 
extend directly to more general situations, however the above 
conditions are exactly the case that we will need
in applications and it is convenient to  restrict our attention to it.
We give the following natural definition: \\

\noindent \bf Definition 5.1. \it We say that the operator 
$T_{n}$ converges on functions $f_{1},...f_{n}$ for some $\lambda$ 
if  the expressions
\[  T^{D_{1},...D_{n}}_{n}(f_{1},...f_{n})(\lambda) \]
converge to a  finite limit as $\min_i D_{i}$ tends to infinity. 
Namely, there exists a number, which we 
denote $T_{n}^{\infty,...\infty}(f_{1},...f_{n})(\lambda)$ such that
for any $\delta >0$ there exists $N_{\delta}$ such that whenever
$\min_i D_{i} \geq N_{\delta}$,  
we have 
\[ \left| T_{n}^{\infty,...\infty}(f_{1},...f_{n})(\lambda) -
 T^{D_{1},...D_{n}}_{n}(f_{1},...f_{n})(\lambda) \right| < \delta. \]

\rm Our first result is the following maximal estimate: \\

\noindent \bf Theorem 5.2. \it Let 
\[ M_{T}(f_{1},...f_{n})(\lambda) = {\rm sup}_{D_{1},...D_{n}}
|T_{n}^{D_{1},...D_{n}}(f_{1},...f_{n})(\lambda)|. \]
Then for any $f_{i} \in L^{p},$ $1\leq p<2,$ we have 
\begin{equation}
\|M_{T}(f_{1},...f_{n})\|_{s_{n}} \leq C_{n}\prod\limits
_{i=1}^{n}\|f_{i}\|_{p}, 
\end{equation}
where $s_{n}^{-1} =nq^{-1}$ and $q^{-1}+p^{-1}=1.$  \\

\begin{proof} \rm The proof uses a well-known device going 
back to Kolmogorov and Seliverstov. Namely, it is sufficient to 
show that for any measurable $D_{i}(\lambda),$ $i=1,...n,$ we have 
\[
\|T_{n}^{D_{1}(\lambda),...D_{n}(\lambda)}(f_{1},...f_{n})\|_{s_{n}}
\leq C_{n} \prod\limits_{i=1}^{n} \|f_{i}\|_{p} \]
with the constant $C_{n}$ independent of $D_{i}(\lambda).$ But  the 
expression $T_{n}^{D_{1}(\lambda),...D_{n}(\lambda)}(f_{1},...f_{n})
(\lambda)$ is obtained from $T_{n}(f_{1},...f_{n})(\lambda)$ by 
replacing the kernels $k_{i}(\lambda,x)$ with the kernels 
\[ \tilde{k}_{i}(\lambda,x)=  k_{i}(\lambda,x) \chi_{R^{+}}(
D_{i}(\lambda)-x_{i}). \]
By Theorem 3.1, these kernels satisfy the estimates (15) for
all $p,$ $q$ such that $p<2$ and $q^{-1}=1-p^{-1},$ with constants
in the norm bounds that do not depend on $D_{i}(\lambda).$ 
Therefore Theorem 4.1 is applicable and directly leads to (33). 
\end{proof} \rm

As one may expect, the maximal estimate (33) implies a.e. convergence. \\

\noindent \bf Theorem 5.3. \it The operator $T_{n}(f_{1},...f_{n})$ 
converges for almost every $\lambda$ on any functions $f_{i} \in L^{p},$ 
$1\leq p <2.$ \\

\begin{proof} \rm Suppose, on the contrary, that for some 
$\epsilon>0$ we have a set 
$S_{\epsilon}$ of positive Lebesgue measure such that for any $N$ there
exist $D_{ij}>N,$ $i=1,...n,$ $j=1,2$ such that 
\[ |T_{n}^{D_{11},...D_{n1}}(f_{1},...f_{n})(\lambda)-
T_{n}^{D_{12},...D_{n2}}(f_{1},...f_{n})(\lambda)|>\epsilon \]
for every $\lambda \in S_{\epsilon}.$ Let us denote 
$f_{i,N}(x_{i})= f_{i}(x_{i})\chi_{R^{+}}(N-x_{i}).$ For any $D_{i}>N,$
we have 
\[ |T_{n}^{D_{1},...D_{n}}(f_{1},...f_{n})(\lambda)-
T_{n}^{D_{1},...D_{n}}(f_{1,N},...f_{n,N})(\lambda)| \leq 
\sum\limits_{j=1}^{n} |T_{n}^{D_{1},...D_{n}}(g_{1j},...g_{nj})(\lambda)|,
\]
where $g_{ij}=f_{i,N}$ if $i < j,$ $g_{ij}=f_{i}$ if $i>j,$ 
and $g_{ij}=f_{i}-f_{i,N}$ if
$i=j.$ In other words, we expand the difference into a telescopic sum.
To estimate each term in the sum, we can apply Theorem 5.2. We get
\begin{eqnarray*} 
&&
{\rm sup}_{D_{1},...D_{n}}\|T_{n}^{D_{1},...D_{n}}(f_{1},...f_{n})
(\lambda)-T_{n}^{D_{1}, \dots D_{n}}(f_{1,N},...f_{n,N})
(\lambda)\|_{s_{n}}^{s_{n}} 
\\
&& 
\qquad\qquad\qquad\le
\sum\limits_{j=1}^{n}\|
T_{n}^{D_{1}, \dots D_{n}}(g_{1j}, \dots g_{nj})(\lambda)\|_{s_{n}}^{s_{n}} 
\\
&&
\qquad\qquad\qquad\le
C_{n}\sum\limits_{j=1}^{n}\left( \prod\limits_{i=1}^{n}\|g_{ij}\|_{p}
\right)^{s_{n}}. 
\end{eqnarray*}
Clearly, the right-hand side goes to zero as $N \rightarrow \infty,$ 
since every product contains the norm of $f_{i}-f_{i,N}$ for some $i.$
On the other hand, the left-hand side is by assumption bounded from below
by $(\frac{\epsilon}{2})^{s_{n}}|S_{\epsilon}|$ for every $N.$ This gives
a contradiction. 
\end{proof} \rm

A straightforward adjustment of the above argument allows
us to pass for almost every $\lambda$ to the 
infinite limit in any order (e.g.\ first in $D_{n},$ then 
in $D_{n-1}$ and so on). 

In order to prove a final Lemma that we will need in the iteration process,
we need to consider a smaller class of multilinear operators than we did
before. These are exactly the operators that appear in the process 
of solving equation (11) by iteration.  \\

\noindent \bf Definition 5.4. \it We say that the multilinear transform 
$T_{n}$ belongs to the class $M_{n},$ if 
\begin{equation}
T_{n}(f_{1},...f_{n})(\lambda)= \int\limits_{0}^{\infty} ...
\int\limits_{0}^{\infty}dx_{1}...dx_{n} \prod\limits_{j=1}^{n}k_{j}
(\lambda, x_{j})f_{j}(x_{j})\prod\limits_{j=2}^{n} \chi_{R^{+}}
(x_{j}-x_{\sigma(j)}), 
\end{equation}
where the kernels $k_{j}(\lambda,x)$ satisfy (15), and the function $\sigma$ 
takes values in $1,...n$ and satisfies $\sigma (j)<j$ for every $j$. \\


\rm
Recall the notation $f_{i,N}(x) = f_i(x)\chi_{R^+}(N-x)$. \\

\noindent \bf Lemma 5.5. \it 
Suppose that for some $\epsilon>0$ and $1\le p < 2$,
for each index $1\le i\le n$,
a function $f_i$ satisfies 
$(1+x)^{ \epsilon}f_{i}(x)\in L^p$.
Then for every $\delta<\epsilon$,
for almost every $\lambda\in R$,
there exists $C(\lambda,n)<\infty$ such that
for every nonnegative $N\in R$,
\[
\left| T_{n}(f_1 - f_{1,N},f_{2} \dots f_{n})(\lambda)\right| 
\leq C(\lambda,n)(1+N)^{-n\delta}.
\]

\begin{proof} \rm 
The structure of operators of class $M_n$ is such that
the integral defining $T_n$ extends only over those
$x$ satisfying $x_j\ge x_1$ for every $j$, 
because of the requirement that $\sigma(j)<j$ for every
$j>1$.
Since $f_1- f_{1,N}$ is supported where $x>N$,
we therefore have for every $N,\lambda$
\[
T_{n}(f_1 - f_{1,N},f_{2} \dots f_{n})(\lambda)
=
T_{n}(f_1 - f_{1,N},f_{2} - f_{2,N} \dots f_{n}
- f_{n,N})(\lambda).
\]
Fix $f_1,\dots f_n$.
Since $\|f_j - f_{j,N}\|_p \le C\cdot(1+N)^{-\epsilon}$
by hypothesis, 
Theorem 5.2 now yields
\[
2^{nr\epsilon} \| \sup_{N\ge 2^r}
|T_n(f_1-f_{1,N},f_2,\dots f_n)|(\lambda)
\|_{L^{s_n}(d\lambda)}
\le C < \infty
\]
for every nonnegative integer $r$, with $C$ independent of $r$.
Consequently
\[
\sum_{r=0}^\infty 2^{nr\delta}
\sup_{N\ge 2^r}
|T_n(f_1-f_{1,N},f_2,\dots f_n)|(\lambda)
< \infty
\]
for almost every $\lambda$. 
To obtain the conclusion of the lemma, given any $N\ge 1$,
choose $r$ so that $2^r\le N<2^{r+1}$ and apply this inequality.
\end{proof} \rm

\vspace{0.2cm}

\begin{center}
 \large \bf 6. Conclusion of the proof of main results.
\end{center}

\vspace{0.2cm}

With the general machinery built up in Sections 4 and 5 in hand, 
we are now in a position to complete the proofs of our main results. 
 First we prove Thereom 1.3. \\

\begin{proof} \rm We recall that it
 suffices to show that under the assumptions of the theorem for a.e.
 $\lambda \in S$ we can find a function 
$q(x,\lambda)$ such that $q(x,\lambda) \in C^{1},$ $q(x,\lambda) 
\rightarrow 0,$ and the condition (9) holds:
\begin{eqnarray*}
q'(x,\lambda) & + & 
\frac{i}{2\Im(\theta
 \overline{\theta})}V(x)\overline{\theta}^{2}(x,\lambda)
\exp (2ip(x,\lambda)) 
\\
&&\qquad + \frac{i}{2\Im(\theta \overline{\theta})}V(x) 
\theta^{2}(x,\lambda)\exp(-2ip(x,\lambda))q^{2}(x,\lambda) 
\qquad \in L^{1}. 
\end{eqnarray*}
By assumption, the kernels 
\[ k_{1}(\lambda,x) =\overline{\theta}^{2}(x,\lambda)\exp (2ip(x,\lambda)) 
\]
and $k_{2}(\lambda, x) = \overline{k_{1}(\lambda,x)}$ satisfy 
$L^{2}-L^{2}$ and trivial $L^{1}-L^{\infty}$ estimates.
Therefore by interpolation, the corresponding operators
map $L^p$ to $L^{q}$, where $q^{-1} =  1 - p^{-1}$
and hence $q>p$, for every $1<p<2$. Therefore 
the theory developed in Sections 4 and 5 applies to the multilinear 
operators $T_{n}$ from classes $M_{n}$ composed from the kernels 
$k_{1}(\lambda, x),$ $k_{2}(\lambda, x)$. 

We will construct the function $q(x,\lambda)$ by iteration. 
Let 
\[ q^{(0)}(x,\lambda) = -\frac{i}{2\Im(\theta \overline{\theta})}
\int\limits_{x}^{\infty}\overline{\theta}^{2}(t,\lambda)\exp
 (2ip(t,\lambda))V(t)\,dt, \]
well-defined for almost every $\lambda \in S$ by Theorem 3.1. Given 
$q^{(n-1)}(x,\lambda),$ we define 
\begin{equation}
q^{(n)}(x,\lambda)=  -\frac{i}{2\Im(\theta \overline{\theta})}
\int\limits_{x}^{\infty}\overline{\theta}^{2}\exp (2ip)V(t)\,dt
 - \frac{i}{2\Im(\theta \overline{\theta})}\int\limits_{x}^{\infty}
 \theta^{2}\exp(-2ip)V(t)(q^{(n-1)})^{2}\,dt.
\end{equation}

Define
\[
V_x(y) = V(y)\chi_{R^+}(y-x).
\]
We need the following \\

\noindent \bf Lemma 6.1. \it $q^{(n)}(x,\lambda)$ is equal to a sum
of multilinear transforms of classes $M_{j}$  (composed from the 
kernels $k_{1}(\lambda,x)$ and $k_{2}(\lambda,x)$ and with all arguments
equal to $V_x$) and is defined for almost every $\lambda \in
S.$ Moreover, 
\[ q^{(n)}(x,\lambda) - q^{(n-1)}(x,\lambda) = \sum\limits_{i} 
T_{j_i}(V_x,V,\dots V), 
\]
where the sum is taken over finitely many orders $j_i$,
each of which satisfies $2n-1 \leq j_{i} \leq 2^{n+1}-1.$ 

\begin{proof} \rm We use induction. For $n=0,$ all statements
 are obvious (defining $q^{(-1)}(x,\lambda)$ to be $0$). 
Suppose they are also
true for $m \leq n-1.$ The fact that $q^{(n)}(x,\lambda)$ is a sum 
of multilinear transforms of classes $M_{j}$ with some 
$j$ follows immediately from the induction hypothesis and 
formula (35). The fact that $q^{(n)}$ is well-defined for a.e.\ 
$\lambda \in S$ is then a consequence of Theorem 5.3. 
Note that 
\begin{eqnarray*} 
\lefteqn{
q^{(n)}(x,\lambda)- q^{(n-1)}(x,\lambda) 
}
\\
& = & 
\frac{-i}{2\Im(\theta \overline{\theta})}
\int\limits_{x}^{\infty}\theta^{2}\exp(-2ip)V
(q^{(n-1)}(t,\lambda)-q^{(n-2)}(t,\lambda))
(q^{(n-1)}(t,\lambda)+q^{(n-2)}(t,\lambda))\,dt. 
\end{eqnarray*}
By the induction hypothesis, every term on the right-hand side is a
multilinear transform of the class $M_{j}$, where the order $j$ is
no less then $(2n-3)+2=2n-1$ and no higher than $(2^{n}-1)2+1= 2^{n+1}-1.$ 
\end{proof} 

\noindent \rm Now we can complete the proof of 
Theorem 1.3. If $x^{\epsilon}V(x) \in L^{p},$ 
$p\leq 2,$ Lemma 5.5 implies that for any multilinear operator of order 
$j_{i}$
\[ \left| T_{j_{i}}^{\infty,\dots\infty}(V_x,...V)(\lambda)
 \right| \leq C(\lambda)(1+x)^{-\delta j_{i}} \]
for almost every $\lambda \in S$ and any $\delta<\epsilon.$ 
Notice that  
\begin{eqnarray*} 
(q^{(n)})'(x,\lambda) & + &
\frac{i}{2\Im(\theta \overline{\theta})}
V(x)\overline{\theta}^{2}(x,\lambda)\exp (2ip(x,\lambda)) 
\\
&& \qquad + 
\frac{i}{2\Im(\theta \overline{\theta})}V(x) \theta^{2}(x,\lambda)
\exp(-2ip(x,\lambda))(q^{(n)})^{2}(x,\lambda) 
\\
& = & 
\frac{i}{2\Im(\theta \overline{\theta})}V(x) \theta^{2}(x,\lambda)
\exp(-2ip(x,\lambda))\left((q^{(n)})^{2}(x,\lambda)-
(q^{(n-1)})^{2}(x,\lambda)\right). 
\end{eqnarray*}
Pick $n$ so that 
$(2n-1)\delta>\frac{1}{2}.$ Then by Lemma 6.1  we find that
the expression on the right hand side is absolutely integrable
for almost every $\lambda$,
since it is the product of some $L^{p}$ 
function $V(x)$ with
a factor which for almost every $\lambda$ is
$O((1+x)^{-\frac{1}{2}-\eta})$
for some  $\eta>0$, and hence belongs to $L^q$ where $q^{-1}=1-p^{-1}$.
 
Hence $q^{(n)}(x,\lambda)$ satisfies the condition (9) and therefore  we
 can take $q(x,\lambda)=q^{(n)}(x,\lambda)$ for a.e. $\lambda \in S.$ 
Then the first claim of Theorem 1.3 follows from Theorem 2.2 and Lemma 2.1.
 As a set of full measure in $S$ it suffices to take the set 
where all multilinear transforms $T_{j}(V,...V)$ 
(composed from $k_{1}(\lambda,x),$ $k_{2}(\lambda,x)$) of order not larger 
than 
$2^{n+1}-1$ converge. It 
remains to prove the formula (2) for the asymptotic behavior of the 
eigenfunctions. Notice that  the asymptotic behavior stated in Theorem 2.2 
differs from the asymptotic behavior we need to show to 
prove Theorem 1.3 only by the presence of an additional multiplier 
\begin{equation}
 \exp \left( \frac{i}{2\Im(\theta \overline{\theta})} 
\int\limits_{0}^{x}(1-|q|^{2})^{-1}V\Re(\theta^{2}q\exp(-2ip))\,dt \right)
\end{equation}
in the asymptotic formula for the solutions in Theorem 2.2. 
But note that the limit of the integral 
\[  \int\limits_{0}^{N}(1-|q(t,\lambda)|^{2})^{-1}V(t)
\Re(\theta^{2}(t,\lambda)q(t,\lambda)\exp(-2ip(t,\lambda)))\,dt \]
as $N \rightarrow \infty$ exists for a.e.\ $\lambda.$ Indeed, we can expand 
$(1-|q|^{2})^{-1}$ into absolutely convergent series in $|q^{2}|.$ Then the
 whole expression becomes represented as a sum of multilinear transforms 
$T_{j}^{N,\infty,...\infty}(V_{x},...V)(\lambda).$ 
Starting from some $l,$ the integrand will become 
absolutely integrable over the whole axis for a.e.\ $\lambda$ by Lemma 5.5; 
in the remaining finite sum every term is convergent by Theorem 5.3. 
Therefore, for a.e.\ $\lambda,$ the 
expression (36) can be written as $C(\lambda)(1+o(1))$ and hence can be 
omitted in the asymptotic expression for the solutions. This completes the
proof. 
\end{proof} \rm 

We now prove Theorems 1.1 and 1.2. \\

\begin{proof} \rm 
It remains only to verify the $L^{2}-L^{2}$ bounds
 for the corresponding operator (4). It is convenient to choose 
$\theta(x,\lambda)=\exp (i\sqrt{\lambda}x)$ in the free case and 
$\theta (x,\lambda)$ Bloch functions in the 
periodic case. The corresponding $L^{2}-L^{2}$ bounds were already 
shown in 
\cite{Kis1}. For the sake of completeness and since the argument is not 
very long, we provide here a sketch of the proof for the free case. For 
the periodic case, the proof is analogous 
given standard information on the properties of the Bloch functions, see 
\cite{Kis1} for details. 

Without loss of generality, we can restrict our attention to some interval 
$(a,b),$ $0<a<b<\infty.$ Consider $\phi(\lambda) \in C_{0}^{\infty}
(0,\infty),$ 
$1\geq \phi(\lambda) \geq 0$ such that $\phi(\lambda) =1$ when 
$\lambda 
\in (a,b).$ Clearly it suffices to show $L^{2}-L^{2}$ bound on functions 
of 
compact support for an operator
\[ Kf(\lambda) = \phi(\lambda)\int\limits_{0}^{\infty}\exp 
\left(2i\sqrt{\lambda}x-\frac{i}{\sqrt{\lambda}}\int\limits_{0}^{x}V(t)\,dt
 \right) f(x)\, dx.
\] 
We have 
\[ \|Kf\|_{L^{2}(a,b)}^{2} = \int\limits_{0}^{\infty}
\int\limits_{0}^{\infty}
dxdy f(x)f(y)\int d\lambda \phi(\lambda) \exp 
\left( 2i\sqrt{\lambda}(x-y)-\frac{i}{\sqrt{\lambda}}
\int\limits_{y}^{x}V(t)\,dt \right). \]
Let us denote by $Z(x,y)$ the kernel 
\[ Z(x,y) = \int d \lambda \phi(\lambda) \exp 
\left( 2i\sqrt{\lambda}(x-y)-
\frac{i}{\sqrt{\lambda}}\int\limits_{y}^{x}V(t)\,dt \right). \]
Let us integrate by parts in $\lambda$ in the expression 
for $Z(x,y)$ $N$ times, integrating $\exp (2i\sqrt{\lambda}(x-y))$ and 
differentiating the rest. We obtain 
\begin{eqnarray} 
|Z(x,y)| & \leq & C(\phi, a,b, N) {\rm min}(1,|x-y|^{-N})\left| 
\int\limits_{y}^{x}
V(t)\,dt \right|^{N}  \\ \nonumber & \leq & C(\phi, a,b, N){\rm min}
(1,|x-y|^{-N(1-\frac{1}{q})})\|V\|_{p}. 
\end{eqnarray}
Taking $N$ large enough, for instance such that $N(1-\frac{1}{q})>1,$ 
we see 
that an operator with the kernel $Z(x,y)$ maps $L^{2}$ to $L^{2}$ 
(by Schur's
 test, for example). Therefore also
\[ \|Kf\|^{2}_{L^{2}(a,b)} \leq C(\phi, a,b)\|f\|^{2}_{2}. \]
\end{proof} \rm

We conclude the paper by formulating one simple generalization of 
Theorem 1.3
 (which implies the corresponding generalizations of Theorems 1.1 and 1.2).
 \\

\noindent \bf Theorem 6.2 \it Fix a  potential $V(x).$ Suppose that there
 exists a monotone differentiable function $d(x),$ $d(x)>0,$ $d(x) 
\stackrel{x \rightarrow \infty}{\longrightarrow} 0,$ $d'(x) \leq 0,$ such 
that $V(x)d(x)^{-1} \in L^{p}$ with some $p<2$ and $V(x)d(x)^{N} \in L^{1}$
 for some integer $N.$ Then under the assumptions of Theorem 1.3, all  
conclusions of  Theorem 1.3 hold for perturbation of $H_{U}$ by $V(x).$ \\

\begin{proof} \rm Going through the proofs of Lemma 5.5 and Theorem 
1.3, we substitute $x^{-\epsilon}$ with $d(x).$ Given that $d(x)^{-1}V(x) 
\in L^{p},$ $p<2,$ all proofs go unchanged, and  in the final step instead 
of $V(x)(1+x)^{-\delta (2n-1)
}$ being absolutely integrable we need to check that for some $N,$ 
\[ V(x)d(x)^{N} \in L^{1}. \]
This is exactly what we assumed in the statement of Theorem 6.2. 
\end{proof} \rm

We note one additional particular class of potentials which we may treat 
using Theorem 6.2. Namely, take any $V(x)$ such as in Theorem 1.3, take a 
sequence $\{x_{n}\}$ and insert intervals of arbitrary size $I_{n}$ at 
each 
point $x_{n}.$ Let $\tilde{V}(x)$ be a potential obtained from $V$ by 
adding
 such intervals $I_{n}$ where
 $\tilde{V}$ is zero. Then it is easy to construct a function $d(x)$ with 
the
 properties as in Theorem 6.2. For the details of such construction we 
refer 
to \cite{Kis} (where it was derived in a slightly different context).  


\begin{center}
 \large \bf Appendix: Singular potentials. 
\end{center}

\vspace{0.2cm}

In Appendix, we discuss the preservation of the absolutely continuous 
spectrum
 for potentials with strong local singularities. The proof turns out to be 
almost entirely parallel to the non-singular case, so we mostly sketch the 
arguments with a few exceptions. 
We should note, however, that the result has some interest in it. 
In the
 explicit construction of power decaying potentials such that the 
corresponding Schr\"odinger operators have purely  singular spectrum 
\cite{Rem} one can try to use the possible
 singularity of the potential to get singular spectrum under stronger
 decay 
conditions. The results of the appendix show that such plan does not work 
out,
 at least the fundamental exponent $\frac{1}{2}$ and virtually all results
 we 
have shown before extend to 
the situation where strong local singularities are allowed.

We will consider the potentials from the spaces $l^{p}(L^{1}),$ $1\leq p 
\leq 2,$ with the norm given by 
\[ \|f\|_{l^{p}(L^{1})} = \left( \sum\limits_{n=0}^{\infty} \left( 
\int\limits_{n}^{n+1} |f(x)|\,dx \right)^{p} \right)^{\frac{1}{p}}. \]

The main result we show is \\

\noindent \bf Theorem A.1. \it Suppose that $U$ is continuous periodic 
(in particular, the free case $U=0$ is of course included). Let the 
perturbation $V(x)$ be such that there exist $p\leq 2$ and $\epsilon>0$ 
with $V(x)x^{-\epsilon} \in l^{p}(L^{1}).$ Then 
all conclusions of Theorems 1.2 hold, in particular the absolutely 
continuous 
spectrum $S$ of the operator 
$H_{U}$ is preserved under the perturbation by $V$ and for  almost every
$\lambda \in S$
there exist solutions of the equation (1) with the pure WKB asymptotic 
behavior in the main term (i.e. with asymptotic behavior given by  (2)). \\

\rm As before, Theorem A.1 is a corollary of the following general 
criterion (an analogue of Theorem 1.3): \\

\noindent \bf Theorem A.2. \it  Suppose that the potential $V(x)$ is 
such that there exist $\epsilon>0$ and $p \leq 2$ so that $x^{\epsilon}
V(x)\in l^{p}(L^{1}).$ Assume  that an operator
\[
(Kf)(\lambda) = \chi(S)\int\limits_{0}^{\infty}\theta(x,\lambda)^{2} \exp 
\left(\frac{i}{\Im(\theta\overline{\theta}')}\int\limits_{0}^{x}V(t)|
\theta(t,\lambda)|^{2}\,dt\right)f(x)\,dx 
\]
satisfies the  bound
\[ \|Kf(\lambda)\|_{L^{2}(S)} \leq C \|f\|_{l^{2}(L^{1})} \]
 on functions $f$ of compact support. Then the absolutely continuous  
spectrum
 of $H_{U}$ supported on the set $S$ is 
preserved under perturbation by $V,$ i.e. the set $S$ belongs to the 
essential
 support of the absolutely continuous part of the spectral measure of 
operator
 $H_{U+V}.$
Moreover, for almost every $\lambda \in S,$ there exist solutions 
$\psi_{\lambda}(x),$
 $\overline{\psi}_{\lambda}(x)$ of the equation (1)
with the asymptotic behavior (2).  \\

\rm We will indicate the changes in the proof of Theorem 1.3 which are 
necessary to prove Theorem A.2. First, we need the following substitute 
for Theorem 3.1.  \\

\noindent \bf Lemma A.3. \it Suppose that $1\le p<\infty$,
that $q>p$, and that $K$ is a bounded linear
operator from $l^p(L^1(R))$ to $L^q(R)$. 
Then the maximal function given by 
\[ (M_{K}f)(\lambda) = {\rm sup}_{N}
| K(f \chi_{R^+}(N-x))(\lambda)|
\]
satisfies 
\[ \|M_{K}f\|_{q} \leq C_{q}\|f\|_{l^{p}(L^{1})} \]
for the same pair of exponents $p,q$.
Moreover the integral defining 
$(M_{K}f)(\lambda)$ converges for almost every $\lambda$, 
for any $f \in l^{p}(L^{1})$.

\begin{proof} \rm Given the function $f$ with 
$\|f\|_{l^{p}(L^{1})}=1,$
 we consider the  family of intervals $E(m,l)$
similar to the $L^{p}$ case. Namely, first we consider intervals 
$E(1,1)$
 and $E(1,2),$ such that their union is the whole half-axis, the first 
interval lies to the left of the second and 
\[ \|f\chi_{E(1,1)}\|^{p}_{l^{p}(L^{1})}=
\|f\chi_{E(1,2)}\|^{p}_{l^{p}
(L^{1})}. \]
 We note that in contrast to the $L^{p}$ case, we cannot in general 
say that
  $\|f\chi_{E(1,i)}\|^{p}_{l^{p}(L^{1})}=2^{-1}.$ However, 
it is easy to
 see that 
\[ \|f\chi_{E(1,i)}\|^{p}_{l^{p}(L^{1})} \leq 2^{-1}. \]
If the sets $E(1,1),$ $E(1,2)$ are unions of the integer intervals 
$(l,l+1),$ we have equality. 
Otherwise each set $E(1,i)$ is a union of 
integer intervals and a set $A_{i},$ such that $A_{1} \cup A_{2} =(m,m+1)$ 
for some $m.$ In this case
\[ \left( \int |f\chi_{A_{1}}|\,dx \right)^{p} + 
\left( \int |f\chi_{A_{2}}|\,dx
 \right)^{p} \leq \left( \int\limits_{m}^{m+1}|f(x)|\,dx \right)^{p} \]
(since $p \ge 1$), and hence 
\[ \|f\chi_{E(1,2)}\|^{p}_{l^{p}(L^{1})}+
\|f\chi_{E(1,1)}\|^{p}_{l^{p}
(L^{1})} \leq 1. \]
We decompose each of the intervals $E(1,i)$ further as in 
$L^{p}$ case.
 On the step $m,$ we obtain a family of intervals $E(m,l),$ such that 
$\cup_{l}E(m,l)=R^{+}, $ $E(m,l_{1})$ lies to the left of
$E(m,l_{2})$ if $l_{1} \leq l_{2}.$ We choose these intervals so that 
the total number of the intervals of generation $m$ is $2^{m}$ and 	
$\|f\chi_{E(m,l)}\|^{p}_{l^{p}(L^{1})} \leq 
2^{-m}$ for every $l.$
Such a family is obtained by splitting every interval of given generation
into two equal pieces, with the same arguments as in the first step. 

Let 
\[ (M_{K,m}f)(\lambda)={\rm sup}_{l}| M_{K}(f\chi_{E(m,l)})(\lambda)|. 
\]
Then 
\[ (M_{K}f)(\lambda) \leq \sum\limits_{m=1}^{\infty}(M_{K,m}f)(\lambda). \]
This follows from the construction of the family $E(m,l).$ 
Indeed, modulo a set on which $f$ vanishes almost everywhere,
any interval $[0,N]$ may be decomposed as a disjoint union
of intervals $E(m,l)$, with at most one such interval 
for each generation number $m$. Summing over $m$ and invoking the
triangle inequality leads to the desired majorization for $M_K f$.
Consequently
\[ \|M_{K}f\|_{q} \leq \sum\limits_{m=1}^{\infty}\|M_{K,m}f\|_{q}.  \]

On the other hand, 
\[ \|M_{K,m}f\|_{q}^{q} \leq \sum\limits_{l=1}^{2^{m}} 
\|K(f\chi_{E(m,l)})\|_{q}^{q} \leq \sum\limits_{l=1}^{2^{m}} 
\|f\chi_{E(m,l)}\|^{q}_{l^{p}(L^{1})} \leq 2^{m(1-\frac{q}{p})}. \]
Therefore,
\[ \|M_{K}f\|_{q} \leq \sum\limits_{m=1}^{\infty} 
2^{m\left(\frac{1}{q}-\frac{1}{p}\right)} \leq C_{q}\|f\|_{l^{p}(L^{1})} \]
(we assumed $\|f\|_{l^{p}(L^{1})}=1,$ but the bound extends to all $f$ by 
sublinearity of $M_{K}$). Almost everywhere
convergence follows from the maximal estimate in a standard way.
\end{proof} \rm

Now we prove Theorem A.2. \\

\begin{proof}
\rm The proofs of multilinear transform properties and almost
everywhere convergence estimates go exactly the same way as before.
The family $E(m,l)$ has the same properties as in $L^{p}$ case, in 
particular martingale-type property (two sets either disjoint or one is 
contained in another). The Lemma 4.3 clearly remains valid. 
The function 
\[ f(x) = \left( \sum\limits_{i=1}^{n} |f_{i}(x)|^{p} \right)^{\frac{1}{p}}, \]
used to construct the family $E(m,l)$ in the proof may be replaced 
by 
\[ f(x)=  \sum\limits_{i=1}^{n} |f_{i}(x)|. \]
The only other change we need to make in the proof is to
change throughout $\|\cdot\|_{p}$ to $\|\cdot\|_{l^{p}(L^{1})}.$ 
\end{proof} \rm

To prove Theorem A.1, we need to show the $l^{2}(L^{1})-L^{2}(S)$ bounds. \\

\begin{proof} \rm We will show the proof of the 
needed norm bound only for the case $U=0.$ The general periodic case 
follows from the properties of the Bloch functions in a way parallel 
to the free case. We refer to \cite{Kis1} for necessary 
information and a similar argument. 

Clearly we can restrict our attention to some compact interval 
$I =(a,b),$ $b>a>0.$ It is sufficient to show that the  operator
\[ (Kf)(\lambda)= \phi(\lambda)^{2}\int\limits_{0}^{\infty} 
\exp\left( 2i\sqrt{\lambda}x - \frac{i}{\sqrt{\lambda}}
\int\limits_{0}^{x}V(t)\,dt \right)f(x)\,dx \]
satisifes $l^{2}(L^{1})-L^{2}(I)$ bound on functions of compact support, 
where $\phi \in C_{0}^{\infty}(R^{+})$ and $\phi(\lambda)=1$ for $\lambda 
\in I.$
Note that 
\begin{eqnarray} 
(Kf)(\lambda) & = & \sum\limits_{l=0}^{\infty}\phi(\lambda)
\exp\left(2i\sqrt{\lambda}l-\frac{i}{\sqrt{\lambda}}
\int\limits_{0}^{l}V(t)\,dt\right)
\nonumber
\\
&& \cdot\ 
\phi(\lambda) \int\limits_{0}^{1}\exp\left(2i\sqrt{\lambda}y-
\frac{i}{\sqrt{\lambda}}\int\limits_{l}^{y+l}V(t)\,dt\right)f(y+l)\,dy. 
\end{eqnarray}
We can write $Kf(\lambda)$ as follows:
\[ Kf(\lambda) = \sum\limits_{l=0}^{\infty}\phi(\lambda)
\exp\left(2i\sqrt{\lambda}l-\frac{i}{2\sqrt{\lambda}}
\int\limits_{0}^{l}V(t)\,dt\right)
 f(\lambda,l), \]
where the expression $f(\lambda,l)$ has the following property: for every
$m \geq 0,$ 
\[ |\partial_{\lambda}^{m}f(\lambda,l)| \leq C_{m}(I,\phi){\rm sup}_{n} 
\left(\int\limits_{n}^{n+1}|V(x)|\,dx\right)^{m}f_{l} \]
where $f_{l} = \int_{l}^{l+1}|f(x)|\,dx.$ This property 
of $f(\lambda,l)$ is evident from (38).

We compute
\[ \|Kf\|^{2}_{L^{2}(I)} \leq \sum\limits_{l=0}^{\infty} 
\sum\limits_{m=0}^{\infty} \int d\lambda \phi^{2}(\lambda)
\exp\left(2i\sqrt{\lambda}(l-m)-\frac{i}{2\sqrt{\lambda}}
\int\limits_{m}^{l}V(t)\,dt\right)f(\lambda,l)\overline{f}(\lambda,m). \]
Let us integrate by parts, differentiating 
$f(\lambda,l)\overline{f}(\lambda,m)$ and integrating the rest. 
By virtually the same computation as one which led us to (37) we obtain
\[
  \left| \int d\lambda \phi(\lambda)\exp
\left(2i\sqrt{\lambda}(l-m)-\frac{1}{2\sqrt{\lambda}}
\int\limits_{m}^{l}V(t)\,
dt\right) \right| \leq C_{r}{\rm min}(1,|l-m|)^{-\frac{r}{p}} 
\]
for every positive integer $r.$  Therefore, taking into account the 
properties
 of $f(\lambda,n),$ we obtain
\[ \|Kf\|^{2}_{L^{2}(I)} \leq C_{r}\sum\limits_{l=0}^{\infty} 
\sum\limits_{m=0}^{\infty} {\rm min}(1,|m-l|)^{-\frac{r}{p}}f_{m}f_{l}. \]
Taking $r$ large enough ($\frac{r}{p}>1$ will do) we see that the operator
$K$ satisfies the required bound. 
\end{proof}

\newpage

\begin{center}
\large \bf Acknowledgment.
\end{center}
A.K. thanks B.~Simon for 
constant support and encouragement and expresses his gratitude to 
S.~Molchanov and C.~Remling for stimulating discussions. 
A.K.'s work at the MSRI has been supported in part by NSF grant DMS 902140. 
M.C.'s work has been supported in part by NSF grant DMS96-23007.

\end{document}